\documentclass[11pt,a4paper]{article}

\usepackage{amsmath}
\usepackage{amsthm}
\usepackage{amsfonts}
\usepackage{epsfig}
\usepackage{cmap}

\usepackage[english]{babel}

\parskip\smallskipamount
\newtheorem{theorem}{Theorem}
\newtheorem{lemma}{Lemma}
\theoremstyle{remark}
\newtheorem{rem}{Remark}

\renewcommand{\phi}{\varphi}
\renewcommand{\epsilon}{\varepsilon}

\parskip \smallskipamount

\title{Doubly Randomised Protocols for a Random Multiple Access Channel with ``Success-Nonsuccess'' Feedback}

\author{Sergey Foss,\thanks{Postal adress: School of Mathematical and Computer Sciences, Heriot-Watt University,
EH14 4AS, Edinburgh, United Kingdom. Email address:
s.foss@hw.ac.uk. Author thanks grant 0770/GF3 of Kazakhstan
Ministry of Education and Science for the support.} {\em
\normalsize Heriot-Watt University, Sobolev Institute of Mathematics and Novosibirsk State University}
\and Bruce Hajek,
\thanks{Postal address: CSL, College of Engineering, University
Illinois at Urbana- Champaign, Urbana, IL 61801ñ2307 USA. Email
address: b-hajek@uiuc.edu} {\em \normalsize University of
Urbana-Champaign}, \and Andrey Turlikov,\thanks{Postal address:
Dept of Information Systems and Data Protection, St.-Petersburg
University of Aerospace Instrumentation, 67, B. Morskaya st.
St.-Petersburg, 190000, Russia. Email address: turlikov@vu.spb.ru.
Author thanks grant 0770/GF3 of Kazakhstan Ministry of Education
and Science for the support.} {\em \normalsize Sankt-Petersburg
State University of Aerospace Instrumentation}}
\date{\today}

\begin{document}

\maketitle

\begin{quotation}\small
We consider a model of a decentralized multiple access system with
a non-standard binary feedback where the empty and collision situations
cannot be distinguished. We show that, like in the case of a ternary feedback,
for any input
rate $\lambda < e^{-1}$, there exists a ``doubly randomized'' adaptive
transmission protocol which stabilizes the behavior of the system.
We discuss also a number of related problems and formulate some
hypotheses.
\end{quotation}


\section{Introduction}
\label{sec:introduction}

We consider a decentralised multiple access system model with an
infinite number of users, a single transmission channel, and
an adaptive transmission protocol that does not use the individual
history of messages. We consider a class of protocols where the user cannot
observe the individual history of messages and the total number of messages.
With any such a protocol, all users transmit
their messages in time slot $(n,n+1)$
with equal probabilities $p_n$ that depend on the history of
feedback from the transmission channel.

Algorithms with ternary feedback ``Empty-Success-Collision''  were
introduced in \cite{Tsy78(1)} and  \cite{Cap79}.  It is assumed
that the users can observe the channel output and distinguish
among three possible situations: either no transmission (``Empty'')
 or transmission from a single server (``Success'') or a collision
 of messages from two or more users (``Conflict'').
It is known since the 80's (see e.g. \cite{Hajek}, \cite{Mih})
that if the feedback is ternary, then the channel capacity is
$e^{-1}$: if the input rate is below $e^{-1}$, then there is a
stable transmission protocol; and if the input rate is above
$e^{-1}$, then any transmission protocol is unstable. A stable
protocol may be constructed recursively as follows: given
probability $p_n$ in time slot $(n,n+1)$ and a feedback at time
$n+1$, probability $p_{n+1}$ is bigger than $p_n$ if the slot
$(n,n+1)$ was empty, $p_{n+1}=p_n$ if there was a successful
transmission, and $p_{n+1}$ is smaller than $p_n$ if there is a
conflict. Hajek \cite{Hajek} considered a multiplicative
increase/decrease and assumed that the random number, $\xi$, of arrivals
per typical slot has a finite exponential moment, ${\mathbf E} e^{c\xi}<\infty$,
for some $c>0$. Mihajlov \cite{Mih} considered an additive increase/decrease
and assume the second moment ${\mathbf E} \xi^2$ to be finite. Later Foss
\cite{Foss} showed that, without further assumptions on the input,
condition $\lambda ={\mathbf E} \xi < e^{-1}$
is sufficient for existence of an multiplicative stable algorithm, and extended
this result onto
a more general stationary ergodic input.

It is also known, see e.g. \cite{Meh84} and \cite{Par92} that similar results (existence/nonexistence
of a stable protocol if the input rate is below/above $e^{-1}$) hold
for systems with either ``Empty--Nonempty'' binary feedback or
``Conflict--Nonconflict'' feedback.

In this paper we show that the channel capacity is again $e^{-1}$
for a third type of the binary feedback 
which may be called ``Success--Nonsuccess''  feedback.
The problem is interesting from a practical point of view,  because
in order for a receiver to distinguish between ``collision'' and ``no
transmission'', it would have to differentiate between the increased
energy or structure present when there is a collision of two or more
packets, from thermal noise. That can be difficult or impossible for
some receivers to do.

Tsybakov and Beloyarov, \cite{Tsy90(4)r} and \cite{Tsy90(5)r},
introduced and studied a model with the ``success-failure'' feedback,
but with an extra option. There is a selected in advance station
which works as follows. Given a nonsuccess feedback, the station
may  send, in the next time slot, a testing package
to recognize what has happened, either empty slot
or collision. Clearly, if an algorithm uses this option regularly,
it is an algorithm with the ternary feedback which works slower that
the conventional one
(uses two time slots instead of one in the case of nonsuccess).
In \cite{Tsy90(4)r} and \cite{Tsy90(5)r}, the authors introduce
a class of algorithms that send a testing package from time to time
only
and show (numerically) that the lower is the rate of using this
option, the closer to $e^{-1}$ is the throughput.

In this paper, we do not allow thus use of test packets.
Our approach to the problem is to introduce a further (second)  randomization.
We consider
a new class of ``doubly randomised'' protocols and
show that, for any pair of numbers $0 <
\lambda_0 < \lambda_1 < e^{-1}$, there exists a protocol
from the
class
that makes stable a system with input rate $\lambda$,
for any $\lambda \in [\lambda_0,\lambda_1]$.
Then we
formulate two conjectures on stability of other classes of protocols
and, in particular of protocols that do not depend on an actual
value of $\lambda$.
Our stability result is based on a generalized
Foster criterion and the fluid approximation approach,
see e.g. \cite{FK}.

In a recent paper \cite{TuFo} (see also \cite{Mal95}), a stability result has been
obtained for a similar model with ``success-failure'' feedback
where a user may also take into account the arrival time of its message.
It was shown that if the input rate is below $0,317$, then there
exist stable algorithms (called ``algorithms with delayed
intervals''). The results of our paper are stronger in two
directions: we show that a stable protocol exists if the input
rate is below $e^{-1}$ (clearly, $e^{-1}>0.317$) and that there is
no need to use the information about arrival times.

There is an interesting question which seems to be open: assume we know
arrival times of messages. Could this extra
information increase an actual channel capacity?

The paper is organised as follows. Section 2 contains the description of the model and
of the class of transmission protocols under consideration, as well as the statement of
the main result. Its proof is presented in Section 3. Then in Section 4 we introduce
two more classes of protocols and formulate conjectures on their stability.

\section{The Model and the Class of Protocols}
\label{sec:Model_algorithm}

We consider (a variant of) a multi-access system introduced in
\cite{Tsy80(1)}. There is an infinite number of users and a single transmission channel available
to all of them. Users exchange their messages using the channel. Time is slotted and all message
lengths are assumed to be equal to the slot length (and equal to one).

The input process of messages $\{ \xi_n\}$ is assumed to be i.i.d., having a
general distribution
with finite mean $\lambda = {\mathbf E} \xi_1$, here $\xi_n$ is a total number of messages arriving within
time slot $[n,n+1)$ (we call it ``time slot $n$'', for short).

The systems operates according to an ``adaptive ALOHA protocol'' that
may
be described as follows.
There is no coordination between the users, and at the beginning of time slot $n$ each message
present in the system is sent to the channel with probability $p_n$, independently of everything
else.   So given that the total number of messages is $N_n$, the number of those sent to the channel,
$B_n$, has conditionally the Binomial distribution $B(N_n,p_n)$ (here $B_n \equiv 0$ if $N_n=0$). Let
$J_n=1$ if $B_n=1$ and $J_n=0$, otherwise. If $J_n=1$, then there is a successful transmission
within time slot $n$. Otherwise there is either an empty slot ($B_n=0$) or a collision of messages
($B_n \ge 2$), so there is no transmission. Then the
following recursion holds:
\begin{equation}\label{key}
N_{n+1}= N_n-J_n+\xi_{n+1}.
\end{equation}

A transmission protocol is determined by sequence $\{p_n\}$. We consider ``decentralised''
protocols: the numbers $N_n$, $n=1,2,\ldots$ are not observable, and only values of past $J_k, k<n$
are known.
We consider protocols where $\{p_n\}$ are defined recursively in the Markovian fashion:
$p_n$ is a (random) number that depends on the history of the system only through $p_{n-1}$ and
$J_{n-1}$. Then a 2-dimensional sequence $(N_n,p_n), \ n=1,2,\ldots$ forms a time-homogeneous
Markov chain.

In the paper, we introduce three classes of decentralised protocols,
prove a stability theorem for the first class and conjecture similar
results for the two others.
To describe these protocols, we introduce additional notation.

Let $N_1\ge 0$ be the initial number of messages in the system and
$S_1\ge 1$ a positive number (which is an "estimator" of unknown
$N_1$). Let further $\beta \in (0,1),$ $C>0$ and $D>0$ be three
positive parameters, and let $\{I_n\}$ be an i.i.d. sequence that
does not depend on the previous r.v.'s, with ${\mathbf P}
(I_{n}=1) = 1-{\mathbf P} (I_{n}=0)=1/2$.

\begin{rem} In what follows, one can assume a sequence of
estimators $\{S_n\}$ to take integer values only, by assuming that
$D$ and $C$ are integer-valued. But this is not needed, in
general.
\end{rem}

The
class ${\cal A}_1$ of algorithms is determined by $\beta
$, $C>0$, $D>0$, $\{J_n\}$ and $\{I_n\}$ as follows. The
transmission probabilities $p_n$ and the numbers $S_n$ are updated
recursively: given $S_n$, we let
\begin{equation*}
p_n=
\begin{cases}
\beta /S_n & \text{if} \ \ I_n=0,\\
1/S_n & \text{if} \ \ I_n=1,
\end{cases}
\end{equation*}
and then
\begin{equation*}
S_{n+1}=
\begin{cases}
S_n + C & \text{if} \ \ J_n=0,\\
S_n + CD & \text{if} \ \ J_n=1 \ \ \text{and} \ \ I_n=0,\\
\max (S_n - CD,1) & \text{if} \ \ J_n=1 \ \ \text{and} \ \ I_n=1.
\end{cases}
\end{equation*}

In words, the reason for second randomisation is as follows. When we get
a successful transmission, we like keep our estimate $S_n$ not
far from the right value of $N_n$ as long as possible. To increase our
chances, we consider a randomised option for probability $p_n$, with taking
two close, but not identical values, $\beta /S_n$ and $1/S_n$.
So if we get success, we increase the $S$-value for the smaller probability
and decrease for the bigger one.

We denote such an algorithm by $A_1(C,D,\beta )\in {\cal A}_1$.
Two other classes of algorithms are defined in Section \ref{conj}.

We can see that, with any algorithm introduced above,
sequence $\{ (N_n,S_n)\}$ forms a time-homogeneous Markov chain.

{\bf Definition 1.} We say that a Markov chain $\{(N_n,S_n)\}$ is {\it
positive recurrent} if there exists a compact set ${\cal K} \in {\cal R}_+^2$
such that
\begin{itemize}
\item
for any pair of initial values $(N_1,S_1)=(N,S)$,
$$
\tau_{N,S} = \min \{n\ge 1 \ : \ (N_n,S_n)\in {\cal K} \} < \infty \quad \mbox{a.s.}
$$
\item further,
$$
\sup_{(N,S)\in {\cal K}} {\mathbf E} \tau_{N,S} < \infty .
$$
\end{itemize}

A Markov chain $\{(N_n,S_n)\}$ is {\it Harris-ergodic} if there is a probability
distribution $\mu$ such that, for any initial value $(N_1,S_1)=(N,S)$, the distributions
of $(N_n,S_n)$ converge to $\mu$ in the total variation,
\begin{equation}\label{TV}
\sup_A |{\mathbf P} ((N_n,S_n)\in A) - \mu (A)|\to 0, \quad n\to\infty
\end{equation}
where the supremum is taken over all Borel sets $A$ in ${\cal R}_+^2$.

It is well-known (see e.g. \cite{MT93})
that a positive recurrent Markov chain $\{(N_n,S_n)\}$ is Harris-ergodic
if it is aperiodic and
there exist a positive integer $m$, a probability measure $\varphi$, and a positive number $c$
such that
\begin{equation}\label{minorant}
{\mathbf P} ((N_m,S_m)\in \cdot \ | \ (N_1,S_n)=(N,S)) \ge c \varphi (\cdot ),
\end{equation}
for all $(N,S)\in {\cal K}$.

{\bf Definition 2.} We say that  Markov chain $\{(N_n,S_n)\}$ is {\it transient} if there is an initial
value $(N_1,S_1)=(N,S)$ such that $N_n+S_n \to \infty$ a.s., as $n\to\infty$.

{\bf Definition 3.} Algorithm $A$ is {\it stable} if the underlying Markov chain
$(N_n,S_n)$ determined by $A$ is Harris-ergodic,
and  {\it unstable} if the underlying Markov chain is transient.

Here is our main result.

\begin{theorem}\label{th1}
Let $0<\lambda_0<\lambda_1<e^{-1}$ be any two numbers. There exist $C>0$, $\beta_1 \in (0,1)$
such that, for any fixed $\beta \in (\beta_1,1)$, there exists $D_0=D_0(\lambda_0,\lambda_1,\beta )$
such that, for any $D\ge D_0$,
algorithm $A_1(C,D,\beta )$ is stable for any input rate
$\lambda \in [\lambda_0,\lambda_1].$
\end{theorem}

\begin{rem}
If $\lambda > e^{-1}$, then any algorithm with a binary feedback is unstable,
since the same result is known for any algorithm with a ternary feedback.
\end{rem}

\section{Proof of Theorem 1}

We need to introduce a number of auxiliary
functions: for positive numbers
$\beta,\lambda,C,D$ and for
$0\le z < \infty$, let

\begin{eqnarray}
j_1(z,\beta ) &=& \frac{\beta z}{2}e^{-\beta z}, \ \ j_2(z) = \frac{ z}{2}e^{-z};\\
j(z,\beta ) &=& j_1(z,\beta )+j_2(z); \\
a(z,\beta ) &=& \lambda - j(z,\beta ); \\
b(z,\beta ) &=& C (1-j(z,\beta )) + CD(j_1(z,\beta )-j_2(z));\\
r(z,\beta ) &=& a(z,\beta )-z b(z,\beta ).
\end{eqnarray}
Cleary, for any $\beta >0$ and as $z\to\infty$,
$j_1(z,\beta ),j_2(z)$ and $j(z,\beta )$ tend to 0,
$a(z,\beta )\to \lambda$, $b(z,\beta )\to C$, and
$r(z,\beta )\to -\infty$.

We rely on the following auxiliary result.

\begin{lemma}
The functions $j(z,\beta ), a(z,\beta ), b(z,\beta )$ and $r(z,\beta )$ satisfy the
following conditions:\\
for any $0<\lambda_0 < \lambda_1 < e^{-1}$, there exists
$\beta_1\in (0,1)$ 
such that, for any $
\lambda \in [\lambda_0,\lambda_1]$, $\beta \in (\beta_1,1)$
and for any
$C\ge C_1:= \frac{\lambda_1+1}{1-e^{-1}}$ and $D\ge D(C)$
(where $D(C)$ is specified in the proof, see equality \eqref{DC} below),
\begin{itemize}
\item{}
equation $a(z,\beta )=0$ has two roots $0<z_1<z_2<\infty$;
\item{}
equation $b(z,\beta )=0$ has two roots $0<t_1<t_2<\infty$;
\item{}
$0<t_1<z_1<t_2<z_2$;
\item{}
function $r$ is continuous in $z\in [0,\infty ]$,
$r(z,\beta )>0$ for $z\le z_1$ and $r(z,\beta )<0$ for $z\ge t_2$;
therefore $\inf_{0\le z \le z_1} r(z,\beta )>0$,
$\sup_{t_2\le z \le \infty} r(z,\beta ) <0$ and all roots to equation $r(z,\beta )=0$ lie
in the interval $(z_1,t_2)$.
\end{itemize}
\end{lemma}

{\sc Proof.} Introduce a further function
$$
b_1(z,\beta ) = 1-j(z,\beta ) + D(j_1(z,\beta )-j_2(z)).
$$
Then
$b(z,\beta ) = C b_1(z,\beta )$.

We know that $ze^{-z}$ (and also $j_2(z)$) is increasing in $z$
if
$z\in (0,1)$ and decreasing if $z>1$. Then, for $\beta <1$, $j_1(z,\beta )$
is increasing if $z<1/\beta$ and decreasing if $z>1/\beta$.
Further, $m(\beta ):= \min_{1\le z\le 1/\beta} j(z,\beta )$ is
strictly positive and tends to $e^{-1}$ as $\beta\uparrow 1$.

For any $\lambda_1 \in (0,e^{-1})$ and any
$\varepsilon \in (0,e^{-1}-\lambda_1)$, one can choose
$\widehat{\beta}_1<1$ so close to 1 that
$m(\beta )\ge \lambda_1+\varepsilon$, for all $\beta\in
[\widehat{\beta}_1,1)$. Then for any
$\lambda \in (0,\lambda_1]$ and any $\beta\in [\widehat{\beta}_1,1]$,
equation $j(z,\beta )=\lambda$ has two roots, $z_1(\beta,\lambda )$
and $z_2(\beta,\lambda )$, with
$$
z_1(\beta,\lambda )\le z_1(\beta,\lambda_1) < 1 <
z_2(\beta,\lambda_1)\le z_2(\beta,\lambda ).
$$
By continuity of functions under consideration,
for any $\lambda\le \lambda_1$, $\beta z_2(\beta,\lambda )\to
z_2(1,\lambda )$ as $\beta\uparrow 1$, so one can choose
$\beta_1\in [\widehat{\beta}_1,1)$ such that
$$
\inf_{\lambda\in (0,\lambda_1]}
\inf_{\beta\in [\beta_1,1]} \beta z_2(\beta,\lambda )>1.
$$
Then, again for all $\lambda \in (0,\lambda_1]$ and all
$\beta\in [\beta_1,1)$,
\begin{equation}\label{jjs}
j_1(z_1,\beta )-j_2(z_1) < 0 < j_1(z_2,\beta )-j_2(z_2),
\end{equation}
with $z_i=z_i(\beta,\lambda )$, for $i=1,2$.

Now we fix $\beta\in [\beta_1,1)$ and, for given
$0<\lambda_0<\lambda_1$, let
\begin{equation}\label{lowerD}
D_0 = \frac{2}{\inf_{\lambda_0\le \lambda \le
\lambda_1}(j_2(z_1)-j_1(z_1,\beta ))} < \infty.
\end{equation}
Then, for any $D\ge D_0$ and for any $\lambda \in [\lambda_0,\lambda_1]$, we have $b_1(z_1,\beta )<0$ and $b_1(z_2,\beta )>0$
(the latter inequality always holds).

Take any $D\ge D_0$ and let $t_1<t_2$ be the roots to equation
$b_1(z,\beta )=0$ (one can easily show that the latter equation has
exactly two roots). Then, clearly, $0<t_1<z_1<t_2<z_2$.

Further, we may choose $C$ such that all roots of equation
$r(z)=0$ lie in the interval $(z_1,t_2)$. Indeed,
 for $z>z_2\ge 1$ and $\lambda \le \lambda_1$,
\begin{eqnarray*}
r(z,\beta )
&=&
\lambda - j(z,\beta ) - Cz(1-j(z,\beta )) - CDz (j_1(z,\beta )-j_2(z)) \\
&\le &
\lambda_1 -C(1-e^{-1}) \le -1
\end{eqnarray*}
if
\begin{equation}\label{lowerC}
C\ge C_1 :=\frac{\lambda_1+1}{1-e^{-1}}.
\end{equation}
Further, for $z\le t_1\le 1$ and $\lambda\ge\lambda_0$,
\begin{eqnarray*}
r(z,\beta )
&>&
\lambda - j(z,\beta ) - Cz(1-j(z,\beta )) \\
&\ge &
\lambda -j(t_1,\beta ) - Ct_1\\
&\ge &
\lambda_0-(1+C)t_1
\end{eqnarray*}
since $j(t_1,\beta )\le t_1e^{-t_1}\le t_1$.
The value of $t_1$ is decreasing to $0$
as $D$ tends to infinity. Therefore, one can choose
$D_1=D_1(C)$ such that $(1+C)t_1\le \lambda_0/2$ for
any $D\ge D_1$, and then
let
\begin{equation}\label{DC}
D(C)=\max (D_0,D_1).
\end{equation}

In more detail,
since $ \sup_{t_0\le t\le 1}
\frac{1-j(t,\beta )}{j_2(t)-j_1(t,\beta )} < \infty
$
for any $0<t_0<1$ and
since
$
\lim_{t\downarrow 0} \frac{1-j(t,\beta )}{j_2(t)-j_1(t,\beta )} = \infty,
$
we may choose $D_1=D_1(C)$ so large that
$$
t(D_1) = \max \{t\in (0,1) \ : \
\frac{1-j(t,\beta )}{j_2(t)-j_1(t,\beta )}\ge D_1\}
$$
satisfies inequality
$t(D_1)\le \frac{\lambda_0}{2(C+1)}.$
Therefore, for any
$D\ge D_1$, we have
$t_1<t(D_1)\le \frac{\lambda_0}{2(C+1)}.$

It is left to comment that $r(z,\beta )>0$ for $z\in (t_1,z_1]$ since
$a(z,\beta )\ge 0$ and $b(z,\beta )<0$, and that $r(z,\beta )<0$ for $z\in [t_2,z_2)$
since $a(z,\beta )<0$ and $b(z,\beta )\ge 0$.

The proof of the lemma is complete.

Now we apply the lemma to the proof of the theorem as follows.

We fix $\beta$ and write $j(z)$, $a(z)$, $b(z)$ and $r(z)$ instead of $j(z,\beta )$, $a(z,\beta )$,
$b(z,\beta )$ and $r(z,\beta )$, for short.

{\bf Step 1.}
We introduce fluid limits for the Markov chain under consideration.
Let Markov Chain $(N_k^{m_0},S_k^{s_0})$ start from initial
values $N_0=m, S_0=s$ and assume that $v:=m+s\to\infty$ and that $m/v\to x, s/v\to y$ where $x+y=1$ and $x,y\ge 0$.
Consider a continuous-time Markov process $(N^m(t),S^s(t))$ where
$$
N^m(t) = N_{[tv]}^m, \ S^s(t)=S_{[tv]}^s
$$
and $[z]$ is the integer part of number $z$. Then we consider a family of weak limits of
these processes, as $v\to\infty$. They are indexed by their initial value $(x_0,y_0)$ with
$x_0\ge 0$, $y_0\ge 0$, $x_0+y_0=1$.
By following the standard scheme (see, e.g., \cite{Ryb}, \cite{Dai}, \cite{Bram}), one
can easily show that each such a limit, say $(\widetilde{N}(t),
\widetilde{S}(t))$ is a Lipschitz function with continuous derivatives, and its
derivatives are the functions $a(z)$ and $b(z)$ that were introduced earlier.
In more detail, for any $x_0,y_0\ge 0, x_0+y_0=1$, for any fluid limit
$(\widetilde{N}(t), \widetilde{S}(t)), t\ge 0$ that starts from initial
value $\widetilde{N}(0)=x_0, \widetilde{S}(0)=y_0$, and for any time $t\ge 0$,
if $\widetilde{N}(t)+\widetilde{S}(t)>0$, we let
$z=z(t)=\widetilde{N}(t)/\widetilde{S}(t)\in [0,\infty ]$. Then
the derivatives
are $d\widetilde{N}(t)/dt = a(z)$ and
$d\widetilde{S}(t)/dt = b(z).$

To see that the derivatives are $a(z)$ and $b(z)$ indeed,
we may find the one-step drift. We have

$$
j(m,s) := {\mathbf E} (J_{n} \ | \ N_{n}=m, S_{n}=s) =
\frac{m\beta}{2s} \left(1-\frac{\beta}{s}\right)^{m-1}
+
\frac{m}{2s} \left(1-\frac{1}{s}\right)^{m-1}.
$$
Then
$$
a(m,s) := {\mathbf E} (N_{n+1}-N_{n} \ | \ N_{n}=m,S_{n}=s )
= \lambda - j(m,s).
$$
In the conditions of Theorem \ref{th1}, and for $s>CD$,
$$
b(m,s) := {\mathbf E} (S_{n+1}-S_{n} \ | \ N_{n}=m,S_{n}=s )
= C (1-j(m,s)) + CD
\left(\frac{m\beta}{2s} \left(1-\frac{\beta}{s}\right)^{m-1}
-
\frac{m}{2s} \left(1-\frac{1}{s}\right)^{m-1}\right).
$$

Then, as $m+s\to\infty$, $m/s \to z$,
$$
j(m,s)\to j(z), \ \ \ a(m,s)\to a(z), \ \ \
b(m,s)\to b(z).
$$

{\bf Step 2.} We have to show next that any fluid limit $(\widetilde{N}(t),\widetilde{S}(t))$
(that starts from initial value $(\widetilde{N}(0),\widetilde{S}(0))$ with
$\widetilde{N}(0)+\widetilde{S}(0)=1$) is {\it stable} in the following sense:
for some $\varepsilon \in (0,1)$, there exists finite time $t_{\varepsilon}$
such that $ \widetilde{N}(t_{\varepsilon})+\widetilde{S}(t_{\varepsilon})\le
1-\varepsilon$.
Then, by the general theory (see e.g. \cite{Bram}),
the positive recurrence of the underlying Markov chains follows.

In the positive quadrant ${\cal R}^2\setminus \{ (0,0)\} = \{ (x,y) \ : \ x,y\ge 0, x+y>0 \}$,
introduce a vector field of ``rates'': the rate from point $(x,y)$ is $(a(z),b(z))$ where
$z=x/y$. From the lemma and Step 1, we may deduce that this vector field is ``self-similar'' (rates do not change along any line with tangent
$z$ that starts from the origin). Functions $(a(z),b(z))$ are continuous in $z\in [0,\infty ]$.
Since $a(0)>0, b(0)>0, a(\infty )>0, b(\infty )>0$ and
the functions change their signs in the ``right'' order
$0<t_1<z_1<t_2<z_2$, we have, in particular, $a(z)<0, b(z)<0$ for $z\in (z_1,t_2)$
and $\inf_z (|a(z)|+|b(z)|) >0$.
Each trajectory of a fluid limit (fluid trajectory, for short) is a ``geodesic'' line with respect to this vector field.

From the properties of function $r$ in Lemma 1, it follows that any fluid trajectory moves towards
the cone ${\cal C} := \{ (x,y) \ : \  x/y \in [z_1,t_2] \}$, hits the cone at some time instant and then
never leaves it and drifts towards the origin. To make this observation more rigorous,  we
%
%
introduce a positive
function $R(z)$ with the following properties: $R(0)=1$ and the pair
$(arctan (z),R(z))$, $0\le z \le \infty$ represents a graph (in the polar coordinates) of a smooth
function that splits the positive quadrant into two domains (where
one of the domains is a convex compact neighbourhood of the origin).

%

\begin{figure}[H]
\centering
\includegraphics[width=1.0\columnwidth]{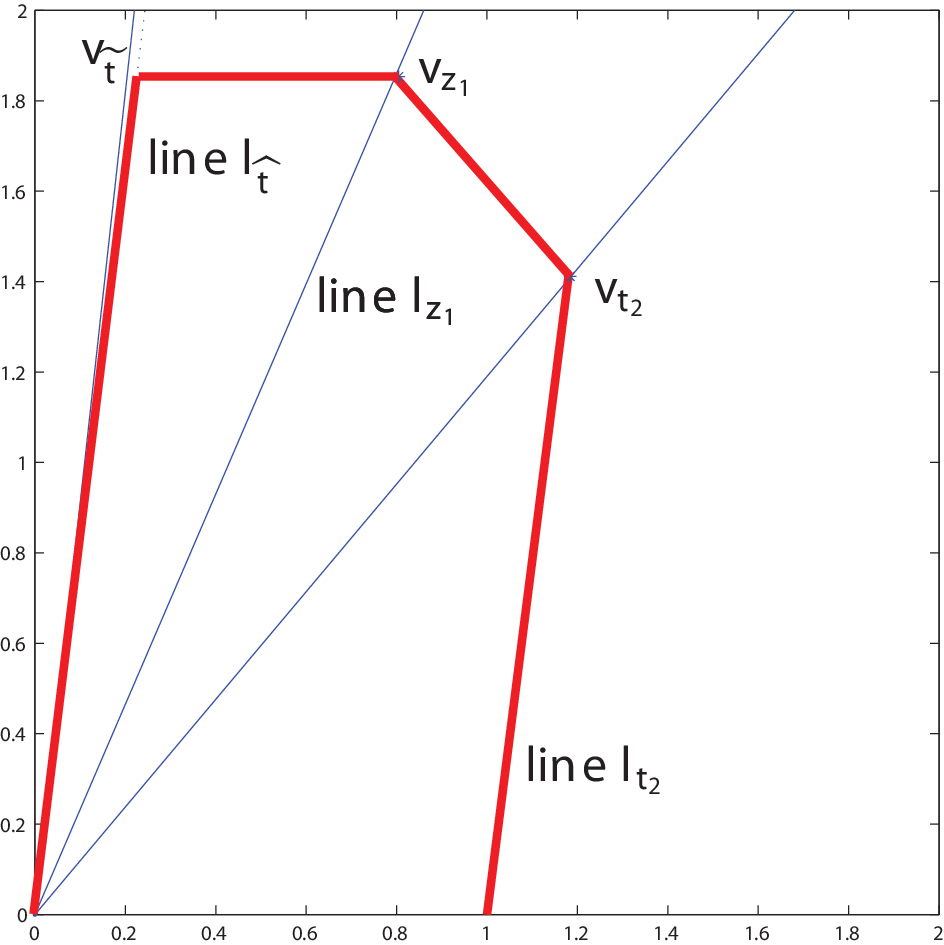}

\caption{Graph L(x,y)=1 in case $\lambda=0.1$ $\beta=0.98$
$C=2.1$ $D=10000$} \label{fig:fig1}
\end{figure}
\bigskip \bigskip

Then, for any constant $c>0$, we draw a line $(arctan (z), cR(z))$, and
introduce a test function $L(x,y)$ by letting $L(x,y)=c$ if
the point $(x,y)$ belongs the line $(arctan (z), cR(z))$.
Further, the line is such that, for $c$ large enough, vector $(a(z),b(z))$ is directed
into the compact domain, and its normal component to the line
is not smaller than a certain positive value.

The construction follows a number of routine steps, therefore we provide
a sketch of the proof and a clarifying figure only.
First we construct a continuous piecewise-linear function. Then we
make it smooth around the points where the line changes
it direction.


We introduce the piecewise linear function as follows.
We start from a fixed point, say $(1,0)$, on the abscissa.
Since $\lambda /C < 1< t_2$, we may choose an angle slightly bigger
than $arctan (\lambda /C)$, say $ \varphi \in (arctan (\lambda /C), \pi /4)$,
and draw a straight line from $(1,0)$ under angle $\varphi$ (measured from the ordinate
line). It intersects the line $l_{t_2}$ with tangent $t_2$ that starts from the origin
(again we measure the tangent
with respect to the ordinate line) at some point, say $v_{t_2}$. From this point, we draw
another straight line at the angle $-\pi /4$ until it crosses the line $l_{z_1}$ with tangent
$z_1$ (that starts from the origin) at some point, say $v_{z_1}$.

Now
we recall that $\lambda /C > t_1$, see the proof of the lemma. Therefore one can take any
 $\varepsilon \in (0, (\lambda /C -t_1)/2)$, let $\widehat{t}= t_1+\varepsilon$ and draw a line
$l_{\widehat{t}}$
with tangent $\widehat{t}$ from the ordinate line
that starts from the origin. Starting from point $v_{z_1}$, we draw a horizontal line in the left direction,
and it intersects with $l_{\widehat{t}}$ at some point, say $v_{\widehat{t}}$.
Finally, starting from $v_{\widehat{t}}$, we draw a straight line under the angle $arctan (\lambda /C - \varepsilon )$
in the left-down direction until it crosses the ordinate line, say at point $v_0$.

Therefore, we have drawn a continuous and convex piece-wise linear line. Then we make it smooth (say
differentiable) by changing in small neighbourhoods of the corners (around points $v_{t_2}, v_{z_1}$ and $v_{\widehat{t}}$),
with keeping it convex.
This completes the construction
of the line $L(x,y)=1$.

Recall that all other lines $L(x,y)=c$
are obtained by the scaling. Then, using routine calculations and the properties of function $r$,
one can show that, for $c$ large enough and for each $z\in [0,\infty ]$,
 the drift vector $(a(z),b(z))$ is directed into the
compact domain $\{ (x,y) \ : \ L(x,y)\le c \}$, and its normal projection is uniformly positive for all $z\in [0,\infty ]$. This implies positive recurrence of the underlying Markov chain.

To conclude that Harris conditions hold we observe that the compact domain contains only
a finite of states (since $N_n$ and $S_n$ are integer-valued), that all these states intercommunicate,
and that the Markov chain is aperiodic since ${\mathbf P} (\xi_1=0)>0$.
This completes the proof of the theorem.

\section{Conjectures}\label{conj}

Here we introduce two more classes of transmission protocols
and conjecture corresponding stability results.

Let ${\cal H}$ be a
class of functions $h:[1,\infty ) \to
[0,\infty )$ such that 
$h(1)=0$, $h(x)\uparrow \infty$ is non-decreasing in $x$, $x-h(x)\uparrow\infty$ is
non-decreasing in $x$,
and $h(x)/x\to 0$, as $x\to\infty$. With each $h\in {\cal H}$, we associate a class
${\cal E}_h$ of positive functions $\varepsilon_h : [1,\infty ) \to (0,1/2]$, such that
$\varepsilon_h(x)\to 0$ and $h(x) \varepsilon^2(x) \to \infty$, as $x\to\infty$.

The two other classes of algorithms differ from the first class in the following.

Algorithms from the second class ${\cal A}_2$ differ from those from
the class ${\cal A}_1$ only in a way the $S$'s are updated: the constant
$CD$ is replaced by a 
function $h\in {\cal H}$. More
precisely, the algorithms are determined by  $\beta$, $C$,
$h(x)$, $\{J_n\}$ and $\{I_n\}$. Given $S_n$, we again
let
\begin{equation*}
p_n=
\begin{cases}
\beta /S_n & \text{if}\ \  I_n=0,\\
1/S_n & \text{if} \ \ I_n=1,
\end{cases}
\end{equation*}
but now define $S_{n+1}$ by
\begin{equation*}
S_{n+1}=
\begin{cases}
S_n + C & \text{if}\ \  J_n=0,\\
S_n + h(S_n) & \text{if} \ \ J_n=1 \ \ \text{and} \ \ I_n=0,\\
\max (S_n - h(S_n),1) & \text{if} \ \ J_n=1 \ \ \text{and} \ \
I_n=1.
\end{cases}
\end{equation*}
For this algorithm, we use notation $A_2( C, h, \beta )\in
{\cal A}_2$.

We modify further the class 2 algorithms by replacing $\beta$ by $1-\varepsilon_h$,
this will form the third class ${\cal A}_3$.
More precisely, the algorithms
are determined by  $C$, $h$, $\varepsilon_h$, $\{J_n\}$ and $\{I_n\}$.
Given $S_n$, we now let
\begin{equation*}
p_n=
\begin{cases}
(1-\varepsilon_h(S_n)) /S_n & \text{if}\ \  I_n=0,\\
1/S_n & \text{if} \ \ I_n=1,
\end{cases}
\end{equation*}
and then define $S_{n+1}$ as for class 2 algorithms:
\begin{equation*}
S_{n+1}=
\begin{cases}
S_n + C & \text{if}\ \  J_n=0,\\
S_n + h_(S_n) & \text{if} \ \ J_n=1 \ \ \text{and} \ \ I_n=0,\\
\max (S_n - h(S_n),1) & \text{if} \ \ J_n=1 \ \ \text{and} \ \
I_n=1.
\end{cases}
\end{equation*}
For this algorithm, we use notation $A_3(C, h,
1-\varepsilon_h)\in {\cal A}_3$.

We can see that again, with any algorithm from the classes
${\cal A}_2$ or ${\cal A}_3$,
a sequence $\{ (N_n,S_n)\}$ forms a time-homogeneous Markov chain.

We believe that the following two statement should be true.

{\bf Conjecture 1.}
Let $0<\lambda_0<e^{-1}$ be any number. There exists $C>0$ and $\beta_2 \in (0,1)$
such that, with any $\beta \in (\beta_2,1)$ and any function $h\in {\cal H}$,
algorithm $A_2(C,h,\beta )$ stabilizes the system, for any input rate
$\lambda < \lambda_0$.\\
If, on the contrary, either $\beta < \beta_2$ or $\lambda > \lambda_0$ ,
then the algorithm $A_2(C,h,\beta )$ is unstable in the system with input rate $\lambda$,
for any
$h\in {\cal H}$.

{\bf Conjecture 2.}
Any algorithm $A_3 (C,h,1-\varepsilon_h)$ from the third class stabilizes the system,
for any input rate $\lambda < e^{-1}$.

\begin{rem}
The conjectures may hold for a broader classes of algorithms if one assumes that,
in the recursion for $S_n$, function $h$ is replaced by two functions, $h_1$ in the second line
and $h_2$ in the third
line.
\end{rem}

\vspace{1cm}


\newpage

\end{document}